\documentclass[12pt]{amsart}
\usepackage{latexsym,amssymb,amsmath}
\usepackage{enumerate}

\newcommand{\Om}{{\Omega}}
\newcommand{\R}{{\mathbb R}}

\newcommand{\Z}{{\mathbb Z}}
\newcommand{\N}{{\mathbb N}}

\newcommand{\C}{{\mathbb C}}
\newcommand{\LL}{{\mathcal L}}

\newtheorem{theorem}{Theorem}[section]

\newtheorem{conjecture}[theorem]{Conjecture}
\newtheorem{lemma}[theorem]{Lemma}
\newtheorem{proposition}[theorem]{Proposition}
\newtheorem{remark}[theorem]{Remark}

\begin{document}
\title[Spectrum is Periodic for $n$-intervals ] {Spectrum is
Periodic for $n$-intervals }

\author[Bose]{Debashish Bose}
\address{Department of Mathematics and
 Statistics, I.I.T. Kanpur, India}
\email{debashishb@wientech.com, madan@iitk.ac.in}

\author[Madan]{Shobha Madan}



\begin{abstract}
In this paper we study spectral sets which are
unions of finitely many intervals in $\R$. We show that any spectrum
associated with such a spectral set $\Omega$ is periodic, with the
period an integral multiple of the measure of $\Omega$. As a consequence
we get a structure theorem for such spectral sets and observe that
the generic case is that of the equal interval case.
\end{abstract}
\maketitle

{\section{Introduction}\label{chap1}}

Let $\Omega$ and $T$ be Lebesgue measurable subsets of $\R^d$ with
finite positive measure. For $\lambda \in \R^d$, let
$$e_{\lambda}(x):=|\Omega|^{-1/2}  e^{2 \pi i \lambda .
x}{\chi}_{\Omega}(x),\,\,\, x\in \R^d.$$

$\Omega$ is said to be a {\bf $spectral$ $set$} if there exists a
subset $\Lambda \subset \R^d$ such that the set of exponential
functions $E_{\Lambda}:=\{e_\lambda:\lambda \in \Lambda\}$ is an
orthonormal basis for the Hilbert space $L^2(\Omega)$. The set
$\Lambda$ is said to be a {\bf $spectrum$} for $\Omega$ and the pair
$(\Omega,\Lambda)$ is called a {\bf $spectral$ $pair$}.
\medskip

$T$ is said to be a {\bf $prototile$} if $T$ tiles $\R^d$ by
translations. In other words, a set $T$ is a prototile if there exists a subset $\mathcal T
\subset \R^d$ such that  $\{T+t: t \in \mathcal T\}$ forms a
partition a.e. of $\R^d$, where $T+t=\{x+t : x \in T\}$. The set
$\mathcal T$ is said to be a {\bf $tiling$ $set$} for $T$ and the
pair $(T,\mathcal T)$ is called a $tiling$ $pair$.
\medskip

The study of relationships between spectral and tiling properties of
sets began with the work of B. Fuglede \cite{Fug}; who while
studying the problem of finding commuting self-adjoint extensions of
the operators $- i \frac{\partial }{\partial x_1}, \cdots, - i
\frac{\partial }{\partial x_n}$ defined on $C_0^\infty(\Omega)$ to a
dense subspace of $L^2(\Omega)$, proved the following result:
\medskip

\begin{theorem}(Fuglede~\cite{Fug})\label{fuglede}
Let $\LL$ be a full rank lattice in $\R^d$ and $\LL^*$ be the dual
lattice. Then $(\Omega,\LL)$ is a tiling pair if and only if
$(\Omega,\LL^*)$ is a spectral pair.
\end{theorem}

He went on to make the following conjecture, which is also known as
the spectral set conjecture:

\begin{conjecture}(Fuglede's conjecture)
{\it A set $\Omega \subset \R^d$ is a spectral set if and only if
$\Omega$ tiles $\R^d$ by translations.}
\end{conjecture}

This led to the study of spectral and tiling properties of sets. In recent years, this
conjecture, in its full generality, has been shown to be false in
both directions if the dimension $d \geq 3$
\cite{T}, \cite{KM1}, \cite{KM2}, \cite{M}, \cite{FR}, \cite{FMM}.
However, interest in the conjecture is alive and the conjecture has
been shown to be true in many cases under additional assumptions.\\

For example, the case where $\Omega$ is assumed to be convex
received a lot of attention recently. It is known that if a convex
body $K$ tiles $\R^d$ by translations then it is necessarily a
symmetric polytope and there is a lattice $\LL$ such that $(K,\LL)$
is a tiling pair \cite{Ven}, \cite{Mac}. Thus the ``tiling implies
spectral'' part of the Fuglede conjecture follows easily from
Fuglede's result. In the converse direction, it has been shown that
a convex set which is spectral has to be symmetric \cite{K2}, and
such sets do not have a point of curvature \cite{IKT1}, \cite{K3},
\cite{IR} (i.e., they are symmetric polytopes). However it is only in
dimension $2$ that the ``spectral implies tiling'' part of the
Fuglede conjecture has been proved \cite{IKT2}.
\medskip

In its full generality Fuglede's conjecture  remains open in
dimensions 1 and 2. In one dimension the conjecture is known to be
related to some interesting number theoretic questions and
conjectures \cite{CM}, \cite{L2}, \cite{LW1}, \cite{Ti}. Even for
the simplest case when $\Omega$ is a finite union of intervals, the
problem is open in both directions and only the $2$-interval case
has been completely resolved by Laba \cite{L1}, where she proved
that the conjecture holds true. In \cite{BCKM} the $3$-interval case
was investigated, where it was shown that for such sets ``tiling
implies spectral'' holds; whereas the ``spectral implies tiling'' part of the
conjecture was proved for this case under some additional
hypothesis.
\medskip

Starting with Fuglede's original work, many results demonstrate that
there  exists a deep relationship between spectra and tiling sets.
For example, when  $I$ is the unit cube  in $\R^d$, then
$(I,\Gamma)$ is a tiling pair if and only if $(I,\Gamma)$ is a
spectral pair. This was first conjectured by Jorgensen and Pedersen
\cite{JP1} who proved it for $d \leq 3$. Subsequently several
authors gave proofs of this result using different techniques
\cite{LRW}, \cite{IP}, \cite{K1}, \cite{Li}. It is worth mentioning
here that tiling by cubes can be very complicated \cite{LS}.
\medskip

In fact there is a dual conjecture due to Jorgensen and Pedersen.

\begin{conjecture}(The dual spectral set conjecture \cite{JP1})
A subset $\Gamma$ of $\R$ is a spectrum for some spectral set
$\Omega$ if and only if it is a tiling set for some prototile $T$.
\end{conjecture}

Approaching the spectral set conjecture by studying the associated
spectra or tiling sets has been very fruitful, specially when these
have some additional structure like periodicity.
\medskip

A set $\Gamma\subset\R^d$ is said to be {\it periodic} if there
exists a full-rank lattice $\LL$ of $\R^d$ such that $\Gamma = \LL +
\{\gamma_1,\dots,\gamma_m\}$, and if, in addition, all coset
differences $\gamma_i-\gamma_j$ are commensurate with the lattice
$\LL$, then $\Gamma$ is said to be {\it rational periodic}.
\medskip

Pedersen \cite{P1} gave a classification of spectral sets which have
a periodic spectrum expressed in terms of complex Hadamard matrices.
On the other hand, Lagarias and Wang \cite{LW2} gave a
characterization of prototiles which tile $\R^d$ by a rational
periodic tiling set in terms of factorization of abelian groups.
Further, they introduced the concept of universal spectrum
\cite{LW2}.
\medskip

A tiling set $\mathcal T$ is said to have a {\it universal spectrum}
$\Lambda_\mathcal T$, if every set $\Omega$ that tiles $R^d$ by
$\mathcal T$ is a spectral set with spectrum $\Lambda_\mathcal T$.
\medskip

Lagarias and Wang \cite{LW2} proved that a large class of tiling
sets $\mathcal T$ have a universal spectrum and then conjectured
that all rational periodic tiling sets have a universal spectrum
which is also rational periodic. This is known as the Universal
Spectrum conjecture. Given a rational periodic tiling set $\mathcal
T$ they gave necessary and sufficient conditions for a rational
periodic spectrum $\Lambda_{\mathcal T}$ to be a universal spectrum
for $\mathcal T$. These developments were instrumental in disproving
the ``tiling implies spectral'' part of Fuglede's conjecture. Later
Farkas, Matolcsi and M\'ora \cite{FMM} proved that the ``tiling
implies spectral'' part of Fuglede's conjecture is equivalent to the
Universal Spectrum conjecture in any dimension.
\medskip

Many results are known concerning the structure of tiling sets
associated with $1$-dimensional prototiles. The fundamental work in
this setting is due to Lagarias and Wang \cite{LW1}, who gave a
complete characterization of the structure of a tiling set $\mathcal
T$ associated with a compactly supported prototile $T$ whose
boundary has measure zero. They proved that in this case $\mathcal
T$ is always rational periodic and the period is an integral
multiple of the measure of $T$. Equipped with this knowledge they
manage to give a characterization of $T$ itself. Further they show
that for every tiling pair $(T,\mathcal T)$ there exists a tiling
pair $(T_1,\mathcal T)$ where $T_1$ is a cluster i.e., a union of
equal intervals, and the problem of finding all possible tiling
pairs $(T,\mathcal T)$ is then related to finding all possible
factorizations of finite cyclic groups. Thus, in essence, the entire
complexity is contained in the equal interval case itself. Later
Kolountzakis and Lagarias extended the periodicity result to all
compactly supported prototiles \cite{KL}.
\medskip

Comparatively much less is known about the structure of spectra
associated with one dimensional spectral sets. All known spectra
associated with one dimensional  spectral sets are rational
periodic. In \cite{JP2} Jorgensen and Pedersen proved that if a
spectral set $\Omega \subset \R$ is a finite union of equal
intervals then it can have finitely many distinct spectra, which are
all periodic. Further, under an additional hypothesis that the set
$\Omega$ is contained in a ``small'' interval, Laba has proved that
the associated spectra for such spectral sets $\Omega$ are rational
periodic \cite{L2}. The general case of spectral sets $\Omega$ which
are unions of finitely many intervals (not necessarily equal) was
studied in \cite{BCKM}. It was shown there that a spectrum $\Lambda$
associated with a spectral set $\Omega$, which is a union of
$n$-intervals has a highly ``arithmetical structure'', namely, if
the spectrum $\Lambda$ contains an arithmetic progression of length
$2n$, then the complete arithmetic progression is contained in it.
\medskip

Our objective in this paper is to study the structure of a spectrum
$\Lambda$ associated with a spectral set $\Omega \subset \R $, when
$\Omega$ is a union of $n$-intervals. We prove that all associated
spectra for such spectral sets are periodic. The essential idea
behind our proof is to show that similar to the case of a tiling set
a finite section of a spectrum essentially determines the complete
spectrum. Theorem \ref{D} and Theorem \ref{orthogonal extension} are
manifestations of this phenomenon and will be central to our proof.
The other key ingredient of the proof is a density result of Landau
for sets of sampling and interpolation (see Theorem \ref{landau}).
In section 2, we state this theorem, explore the geometry of the
zero set of the Fourier transform of a spectral set and prove
Theorem \ref{D} and Theorem \ref{orthogonal extension}.
\medskip

In section 3 we prove our main theorem
\begin{theorem}\label{theorem 1}
Let $\Omega$ be a union of $n$ intervals, $\Omega=\cup_{j=1}^n I_j$,
such that $\left|\Omega\right|=1$. If $\left(\Omega,\Lambda\right)$
is a spectral pair, then $\Lambda$ is a $d$-periodic set with $d\in
\N$. Thus $\Lambda$ has the form $\Lambda=\cup_{j=1}^d
\left\{\lambda_j + d\Z \right\}.$
\end{theorem}

The structure of spectral sets which have a periodic spectrum have
been studied in \cite{P1} and \cite{LW2}. As a consequence of
Theorem \ref{theorem 1} we get a structure theorem for such spectral
sets and observe that the equal interval case is the generic case.

\begin{theorem}\label{theorem 2}
Let $(\Omega, \Lambda)$ be a spectral pair such that $\Omega$ is a
bounded region in $\R$ and $\Lambda$ is $d$-periodic. Then there
exists a disjoint partition  of $[0,1/d)$ into finite number of sets
$E_1,E_2,\dots,E_k$ such that $\Omega = \cup_{j=1}^k (E_j + A_j);\,
A_j \subseteq \Z/d$. Further, each set $\Omega_j:=[0,1/d)+A_j$ is a
spectral set with $\Lambda$ as a spectrum.
\end{theorem}

\section{The geometry of the spectrum}\label{section 2}

Let $(\Omega,\Lambda)$ be a spectral pair. Since spectral properties
of sets are invariant under affine transformations, we will
henceforth assume that $\Omega$ has measure $1$ and that $0\in \Lambda
\subset \Lambda-\Lambda$.
\smallskip

In this paper we will always assume that $\Omega$ is bounded. Then
$\widehat{\chi_\Omega}$, the Fourier Transform of the characteristic
function of $\Omega$, is an entire function.
\smallskip

Let $\Z(\widehat{\chi_\Omega})$ be the zero set of
$\widehat{\chi_\Omega}$ union $\{0\}$ i.e.,
$$\Z(\widehat{\chi_\Omega}) := \{ \xi \in \R:
\widehat{\chi_\Omega}(\xi) = 0 \} \cup \{0\}.$$

If $\lambda, \lambda' \in \Lambda$, then by orthogonality of
$e_\lambda$ and $e_{\lambda'}$ we have $\lambda-\lambda' \in
\Z(\widehat{\chi_\Omega})$. Hence $0\in \Lambda \subset
\Lambda-\Lambda \subset \Z(\widehat{\chi_\Omega})$. Thus the
geometry of the zero set of $\widehat{\chi_\Omega}$ plays a crucial
role in determining the structure of $\Lambda$.
\smallskip

Observe that, as $\widehat{\chi_\Omega}(0)=1$, there exists a
neighborhood around $0$, which does not intersect
$\Z(\widehat{\chi_\Omega})$ except at $0$. Hence, $\Lambda$ is
uniformly discrete. Let $\Lambda_s$ be the set of spectral gaps for a spectrum $\Lambda$
i.e., $$\Lambda_s := \{\lambda_{n+1}-\lambda_n | \lambda_n \in
\Lambda \}.$$

Clearly  $\Lambda_s \subseteq \Lambda - \Lambda \subseteq
\Z(\widehat{\chi_\Omega})$ and $\Lambda_s$ is bounded below. On the
other hand, as a consequence of Landau's density results (see
Theorem \ref{landau} below), we see easily that $\Lambda_s$ is also bounded
above. So, by the analyticity of $\widehat{\chi_\Omega}$ we can conclude that
$\Lambda_s$ is finite. Thus the spectrum can be seen as a
bi-infinite word made up of a finite alphabet, in terms of the
spectral gaps. When $\Omega$ is a union of finite number of
intervals, a much more precise estimate is known for spectral gaps
\cite{Landau}, \cite{IK}, \cite{IP2}.
\smallskip

From now on we will assume that $\Omega$ is a union of a finite
number of intervals. Let $\Omega=\cup_{i=1}^n [a_i,a_i+r_i)$, $\sum_{i=1}^n
r_i=1$. Then,
\smallskip

$$\widehat{\chi_\Omega}(\xi)=\frac{\sum_{i=1}^n [e^{2 \pi i
 (a_i+r_i)\xi} - e^{2 \pi i (a_i)\xi}]}{2\pi i \xi},$$
and $\Z(\widehat{\chi_\Omega})$ is precisely the zero set of the
exponential polynomial given by $$\mathcal
P_\Omega(\xi):=\sum_{i=1}^n (e^{2 \pi i (a_i+r_i)\xi} - e^{2 \pi i
(a_i)\xi}),$$ which is the numerator in the expression of
$\widehat{\chi_\Omega}$. Thus we are naturally led to the study of
exponential polynomials and their zeros.
\smallskip

There is a beautiful result by Turan \cite{Tu}, \cite{nazarov} which
gives size estimates of  exponential polynomials along arithmetic
progressions. This result has the  interesting consequence that if
an arithmetic progression $a,a+d,\dots,a+(2n-1)d$ of length $2n$
occurs in $\Z(\widehat{\chi_\Omega})$ then the complete arithmetic
progression $a+d\Z \subset \Z(\widehat{\chi_\Omega})$. This suggests
that the zero sets of exponential polynomials are highly structured
and we are naturally led to ask the question whether $\Lambda$
inherits this kind of structure?
\smallskip

In the next section we will prove an analog of Turan's Lemma for the
spectrum.

\subsection{Arithmetic Progressions in $\Lambda$}
As we have mentioned before, it was shown in \cite{BCKM} that the
existence of an arithmetic progression of length $2n$ in $\Lambda$
implies that the complete arithmetic progression is in $\Lambda$. Here,
we improve on that result and using Newton's Identities about
symmetric polynomials, give a proof that the occurrence of an
arithmetic progression of length $n+1$ in the spectrum ensures that
the complete arithmetic progression is in the spectrum. Let
$$P(z):= \prod_{i=1}^n (z- \alpha_n) = z^n + S_1
z^{n-1} + S_2 z^{n-2} + \dots + S_n.$$
\medskip

\noindent Let $W_k$ be the sum of k'th power of the roots of $P(z)$,
namely $$W_k:= \alpha_1^k + \alpha_2^k+\dots+\alpha_n^k ;\ \ k=
1,\dots,n. $$
\noindent Then the coefficients $S_i$ and $ W_i$ are
related by the Newton's Identities:
\begin{equation}\label{111}
 W_k + S_1 W_{k-1}+ S_2 W_{k-2} +\dots +
 S_{k-1} W_1 + k S_1 =0;\ \ k=1,\dots,n.
\end{equation}

\noindent Thus $W_1,W_2,\dots,W_n$ uniquely determine the
polynomial $P(z)$.
\medskip

\begin{proposition}\label{C} If $\Z(\widehat{\chi_\Om})$ contains
an arithmetic progression  of length $n+1$ with its first term $0$,
say $0, d, \dots, nd \in \Z(\widehat{\chi_\Om})$ then
\begin{enumerate}[(a)]
\item the whole arithmetic progression $d\Z \subset
\Z(\widehat{\chi_\Om})$,
\item $d\in\Z$, and
\item $\Omega$ d-tiles $\R$.
\end{enumerate}
\end{proposition}
%
%
%
%
\begin{proof}

Note that if $t \in \Z(\widehat{\chi_\Om}) $, then
$$ {\sum}_{j=1}^n  [e^{2\pi i t (a_j + r_j)} - e^{2\pi i t
a_j}] \, = \,0.$$

\noindent The hypothesis says that $\widehat{\chi_\Om}
(ld)=0;\,\, l=1,...,n$, hence

$$ {\sum}_{j=1}^n  [e^{2\pi i ld (a_j + r_j)} - e^{2\pi
i ld a_j}] \, = \,0,\,\, l=1,...,n$$

\noindent We write $\zeta_{2j}=e^{2\pi i d a_j} \, ; \, \,
  \zeta_{2j-1}= e^{2\pi i d (a_j + r_j)};\,\, j=1,...,n$. Then 
the above system of equations can be rewritten as
\begin{equation}\label{333}
\begin{array}{ccccccccccccccc}
\zeta_1+\zeta_3 + \cdots+ \zeta_{2n-1}=\zeta_2+\zeta_4 + \cdots+
 \zeta_{2n}= W_1 \\

\zeta_1^2+\zeta_3^2+\cdots+\zeta_{2n-1}^2=\zeta_2^2+\zeta_4^2
+\cdots+\zeta_{2n}^2 = W_2 \\

\vdots \ \ \ \ \ \ \ \ \ \  \ \ \ \ \ \ \ \ \ \
 \ \ \ \ \ \ \ \ \ \ \vdots \\

\zeta_1^{n}+\zeta_3^{n}+\cdots+\zeta_{2n-1}^{n}= \zeta_2^{n}+
\zeta_4^{n}+\cdots+\zeta_{2n}^{n} = W_n \\
\end{array}
\end{equation}
\medskip

\noindent Let $$P_{1}(z) := \prod_{j=1}^n (z - \zeta_{2j-1})\ \ \
\and\ \ \ \ P_2(z) := \prod_{j=1}^n (z - \zeta_{2j}) $$ Then by
(\ref{111}) and  (\ref{333}) we get $P_1(z) = P_2(z)$.
\medskip

Thus we get a partition of $\zeta_i$'s into n distinct pairs
$(\zeta_i,\zeta_j)$ such that $\zeta_i = \zeta_j; \,\, $  $i\in
1,3,\dots,2n-1$  and $j\in 2,4,\dots,2n$. We can relabel the
$\zeta_{2j}$'s, $j=1,\dots,n$ so that $\zeta_{2j-1}=\zeta_{2j}$. But
then $\zeta_{2j-1}^k = \zeta_{2j}^k, \,\,\forall k \in \Z$ and we
get
\begin{equation}
\widehat{\chi_\Omega}(kd)=\frac{1}{2 \pi i
kd}\sum_{j=1}^n(\zeta_{2j-1}^{k}-\zeta_{2j}^k) = 0;\ \ \forall k \in
\Z \setminus \{0\}.
\end{equation}
Thus $d\Z \subset \Z(\widehat{\chi_\Om})$. Now consider,
\begin{equation}\label{we}
F(x)=\sum_{k\in \Z} \chi_\Omega\left(x+k/d\right),\, x\in [0,1/d)
\end{equation}
Thus $F$ is $\frac{1}{d}$ periodic and integer valued and
\begin{equation}\label{62}
\widehat{F}(ld)= d \sum_{k \in \Z} \int_{0}^{\frac{1}{d}}
\chi_{\Omega}(x+k/d)e^{- 2\pi i ldx} dx =d
\widehat{\chi_\Omega}(ld)=d\, \delta_{l,0}
\end{equation}
Thus $F(t)=d \mbox{ a.e}$. so $ d\in \Z $
and $\Omega$  d-tiles the real line.
\end{proof}

Using Proposition \ref{C}, we now prove the corresponding result for
the spectrum.
\begin{theorem}\label{D}
Let $(\Omega,\Lambda)$ be a spectral pair. If for some $a,d\in
\R,\,$ an arithmetic progression of length $n+1$, say
$a,a+d,...,a+nd \in \Lambda $, then the complete arithmetic
progression $a+d\Z \subseteq \Lambda .$ Further $d \in \Z$ and
$\Omega$ $d$-tiles $\R$.
\end{theorem}
\begin{proof}
Since $  a,a+d,...,a+nd \in \Lambda$, shifting $\Lambda$ by $a$ we
get that $\Lambda_1=\Lambda-a$ is a spectrum for $\Omega$ and
$$0,d,...,nd \in \Lambda_1 \subset \Lambda_1 - \Lambda_1 \subset
\Z(\widehat{\chi_\Om}).$$ Thus surely $d\Z \subset
\Z(\widehat{\chi_\Om}) $ by Proposition \ref{C}.
\medskip

\noindent Now, let $\lambda \in \Lambda_1$. Then by orthogonality,
$$-\lambda,d-\lambda,2d-\lambda,...,nd-\lambda \in
\Z(\widehat{\chi_\Om}).$$
Put
$$\xi_{2j} = e^{-2 \pi i \lambda a_j},\, \xi_{2j-1} = e^{-2 \pi i \lambda(a_j+r_j)}; \,\, j=1,...,n$$
$$\zeta_{2j}=e^{2\pi i d a_j},\, \zeta_{2j-1}=e^{2\pi i
 d(a_j+r_j)}; \,\, j=1,...,n.$$
\\
\noindent Since $\widehat{\chi_\Om}(kd-\lambda)=0$, \, for
$k=0,\dots,n$ we have
\begin{equation}\label{555}
\xi_1 \zeta_1^k - \xi_2 \zeta_2^k + \cdots+\xi_{2n-1} \zeta_{2n-1}^k
- \xi_{2n} \zeta_{2n}^k =
 0 \,\, \mbox{ for }  k=0,...,n.
\end{equation}

\medskip
But the ${\zeta_i}$'s can be partitioned into n disjoint pairs
$(\zeta_i,\zeta_j)$ such that $\zeta_i=\zeta_j$ where $i\in
1,3,\dots,2n-1$  and $j\in 2,4,\dots,2n$. Without loss of
generality, we relabel the $\zeta_{2j}$'s and simultaneously, the corresponding
$\xi_{2j}$'s so that $\zeta_{2j-1}=\zeta_{2j},\, j=1,...,n$. Thus
from (\ref{555}) we get

\smallskip
\begin{equation}\label{666}
\left(\begin{array}{cccc}
1 & 1 & \cdots & 1  \\
\zeta_1 & \zeta_3 & \cdots & \zeta_{2n-1} \\
\vdots & \vdots & \ddots & \vdots  \\
\zeta_1^{n-1} & \zeta_3^{n-1} & \cdots & \zeta_{2n-1}^{n-1}  \\
\end{array}\right)
\left(\begin{array}{c}
\xi_1-\xi_2\\ \xi_3-\xi_4 \\ \vdots \\ \xi_{2n-1} -\xi_{2n}\\
\end{array} \right)=
\left(\begin{array}{c}
0 \\ 0 \\ \vdots \\ 0
\end{array} \right)
\end{equation}
\smallskip

\noindent Now, if $[ \xi_1 -
\xi_2,\xi_3-\xi_4,...,\xi_{2n-1}-\xi_{2n}]^t $  is the trivial
solution, i.e., $\xi_{2j-1}-\xi_{2j}= 0,\,  \forall \,j=1,...,n$
then $\forall \,k \in \Z ,$ we have

$$\widehat{\chi_\Omega}(kd-\lambda) = \frac{1}{2\pi i
(kd-\lambda)} \left[\xi_1 \zeta_1^k - \xi_2\zeta_2^k+
\cdots+\xi_{2n-1}\zeta_{2n-1}^k -\xi_{2n} \zeta_{2n}^k \right]$$
  $$ = \frac{1}{2\pi i (kd-\lambda)} \left[ \zeta_1^k (\xi_1-\xi_2)
  +\cdots+\zeta_{2n-1}^k (\xi_{2n-1}-\xi_{2n}) \right]=0.$$
\smallskip

Thus $d \Z -\lambda \in \Z(\widehat{\chi_\Om})$. If, however,
$[\xi_1-\xi_2,\xi_3-\xi_4,...,\xi_{2n-1}-\xi_{2n}]^t$ is not the
trivial solution, then $\zeta_{2l-1}=\zeta_{2k-1}$ for some $l,k \in
1,...,n ; l\neq k.$

Removing all the redundant variables and writing the remaining
variables as $\eta _{2j+1}^l , \,\, j,l  = 0, 1, \dots, k-1$, we get
a non-singular Vandermonde matrix satisfying
\begin{equation}\label{777}
\left(\begin{array}{cccc}
1 & 1 & \cdots & 1 \\
\eta_1 & \eta_3 & \cdots & \eta_{2k-1} \\
\vdots & \vdots & \ddots & \vdots \\
\eta_1^{k-1} & \eta_3^{k-1} & \cdots & \eta_{2k-1}^{k-1} \\
\end{array}\right)
\left(\begin{array}{c}
\sum_1\\ \sum_3 \\ \vdots \\ \sum_{2k-1}\\
\end{array} \right)=
\left(\begin{array}{c}
0 \\ 0 \\ \vdots \\ 0 \\
\end{array} \right)
\end{equation}
where $$ {\sum}_{k} = \sum_{j:{\zeta}_{2j-1}=\eta_k} {{\xi}_{2j-1} -
{\xi}_{2j}}.$$ Then each of the  $\sum_i =0; \,\, i=1,\dots,k$. But,
then once again $\forall p  \in  \Z,$
$$\widehat{\chi_\Omega}(pd-\lambda)=\frac{1}{2\pi i (pd-\lambda)}
 \left[ \eta_1^p
{\sum}_1 + \eta_3^p{\sum}_3+\cdots +
\eta_{2k-1}^p{\sum}_{2k-1}\right]=0.$$ Thus $d\Z -\lambda  \subseteq
\Z(\widehat{\chi_\Om})$. We already have $d\Z \subseteq
\Z(\widehat{\chi_\Om})$ and now we have seen if $\lambda \in
\Lambda_1$ then $d\Z -\lambda \in \Z(\widehat{\chi_\Om}) $. Thus
$d\Z \subseteq \Lambda_1$, hence $a+d\Z \subset \Lambda$. That $d\in
\Z$ and $\Omega$ $d$-tiles $\R$ follows from Proposition \ref{C}.
\end{proof}

\begin{remark} Theorem \ref{D} is the best possible result in this
direction, as existence of an arithmetic progression of shorter
length in a spectrum does not ensure the complete arithmetic
progression is in the spectrum. For example, consider
$\Omega=[0,1/3]\cup[1,4/3]\cup[2,7/3]$ then
$\Lambda=\{0,1/3,2/3\}+3\Z$ is a spectrum for $\Omega$ which
contains the $3$ term arithmetic progression $0,1/3,2/3$ but clearly
the complete arithmetic progression $\Z/3 \not \subseteq \Lambda$.
\end{remark}

\subsection{Embedding $\Lambda$ in a vector space}
In this section we will investigate the spectrum in a geometric
manner. The setting is again that of a set $\Omega$, which is a
union of finitely many intervals, namely, $\Omega = \cup_1^n [a_j,
a_j+r_j] $. We assume that $\Omega$ is spectral with a spectrum
$\Lambda$. We will embed $\Lambda$ in a vector space and incorporate
the orthogonality of the corresponding set
$E_\Lambda=\{e_\lambda:\lambda \in \Lambda\}$, via a conjugate
linear form.
\medskip

Consider the $2n$-dimensional vector space $\C^n\times\C^n$. We
write its elements as $\underbar{v}=\left(v_1,v_2\right)$ with
$v_1,v_2\in\C^n$. For $\underbar{v},\underbar{w}\in\C^n\times\C^n$
define $$\underbar{v}\odot\underbar{w}:=\langle v_1,w_1\rangle
-\langle v_2,w_2\rangle, $$ where $\langle\cdot,\cdot\rangle$
denotes the usual inner product on $\C^n$. Note that this conjugate
linear form is degenerate, i.e., there exists $\underbar{v} \in
\C^n\times\C^n$, $\underbar{v} \neq 0$ such that
$\underbar{v}\odot\underbar{v}=0$. We call such a vector a {\it
null-vector}. For example, every element of $\mathbb T^n \times
\mathbb T^n$ is a null vector.
\medskip

A subset $S\subseteq\C^n\times\C^n$ is called a set of {\it mutually
null-vectors} if $\forall$ $\underbar{v}, \underbar{w}\in S$ we have
$\underbar{v}\odot\underbar{w}=0$. It is clear from the definition
that elements of a set of mutually null-vectors are themselves
null-vectors.

\begin{lemma}\label{dim V}
Let
$S=\left\{\underbar{v}^1,\underbar{v}^2,\cdots,\underbar{v}^m\right\}$
be a set of mutually null-vectors in $\C^n\times\C^n$. Let $V$ be
the linear subspace spanned by $S$. Then, $V$ is a set of mutually
null-vectors and $dim(V)\leq n$.
\end{lemma}

\begin{proof}
Let $\underbar{v},\underbar{w}\in V$. Since the subspace $V$ is
spanned by $S$, we have $\underbar{v} = \sum_{i=1}^m a_i \
\underbar{v}^i $ and $\underbar{w} = \sum_{j=1}^m b_j \
\underbar{v}^j $. Now, as the set $S$ is a set of mutually
null-vectors $(\underbar{v}^i\odot\underbar{v}^j)=0; \,\, \forall
\,\, i,j=1,\dots,m$ and so, we have
$\underbar{v}\odot\underbar{w}=\sum_{i,j=1}^m a_i \overline{b_j} \
(\underbar{v}^i\odot\underbar{v}^j)=0$. Hence, $V$ is a set of
mutually null-vectors.
\medskip

Let $\underbar{w}^j:=\left(e_j, 0\right),\, j=1,\dots,n$ where
$e_j$'s are the standard basis vectors of $\C^n$. Consider the
subspace $W$ of $\C^n\times\C^n$ spanned by the vectors
$\underbar{w}^j,\, j=1,\dots,n$. Since, these vectors are linearly
independent in $\C^n\times\C^n$, $dim(W)=n$. Further, note that for
$\underbar{w} \in W, \, \underbar{w} \neq 0$ we have
$\underbar{w}\odot\underbar{w} > 0$. Thus $W \cap
V=\left\{0\right\}$ and hence $dim(V) \leq n$.
\end{proof}
\medskip

Suppose $\Omega=\cup_{j=1}^n \left[a_j, a_j+r_j\right)$ is a union
of $n$ disjoint intervals with
$a_1=0<a_1+r_1<a_2<a_2+r_2<\cdots<a_n<a_n+r_n$ and $\sum_1^n r_j=1$.
\medskip

We define a map $\varphi_{_\Omega}$ from $\R$ to $\mathbb
T^n\times\mathbb T^n\subseteq\C^n\times\C^n$ by
$$ x\rightarrow \varphi_{_\Omega}(x)=\left(\varphi_1(x),
\varphi_2(x)\right),$$ where $$ \varphi_1(x)=\left(e^{2\pi i
(a_1+r_1) x}, e^{2\pi i (a_2+r_2) x}, \dots, e^{2\pi i (a_n+r_n)
x}\right)$$ and $$ \varphi_2(x)=\left(1, e^{2\pi i a_2 x}, \dots,
e^{2\pi i a_n x}\right).$$

The following lemma, which is immediate from the definitions, makes
clear the connection between a spectral pair $(\Omega,\Lambda)$ and
the image of $\Lambda$ under the map $\varphi_{_\Omega}$.

\begin{lemma}\label{equivalence}
Let $\Omega$ be a union of $n$ intervals, as above, and suppose
$\Gamma \subseteq \R$ . Then the set of exponentials
$E_{\Gamma}=\{e_\gamma:\gamma \in \Gamma\}$ is an orthogonal set in
$L^2(\Omega)$ if and only if
$\varphi_{_\Omega}(\Gamma):=\left\{\varphi_{_\Omega}(\gamma):\gamma\in
\Gamma \right\}$ is a set of mutually null-vectors.
\end{lemma}

Thus, if $(\Omega, \Lambda)$ is a spectral pair,
$\varphi_{_\Omega}(\Lambda)$ is a set of mutually null-vectors. What
about the converse? We will now try to find some criterion to decide
whether a given pair $(\Omega,\Lambda)$ is a spectral pair.
\smallskip

First, observe that from Lemma \ref{dim V}, we already know that if
$(\Omega,\Lambda)$ is a spectral pair then the vector space
$V_\Omega(\Lambda) :=
span\left\{\varphi_{_\Omega}(\lambda):\lambda\in\Lambda\right\}$ has
dimension at most $n$. We will now show that $\Lambda$ has a ``local
finiteness property'', in the sense that there exists a finite
subset $\mathcal{B}$ of $\Lambda$, $\# \mathcal{B} \leq n$, such
that $\Lambda$ gets uniquely determined by $\mathcal{B}$.

\begin{lemma}\label{local finite}
Let $(\Omega,\Lambda)$ be a spectral pair and
$\mathcal{B}=\left\{y_1,\dots,y_m\right\}\subseteq\Lambda$ be such
that $\varphi_{_\Omega}(\mathcal
B):=\{\varphi_{_\Omega}(y_1),\dots,\varphi_{_\Omega}(y_m)\}$ forms a
basis of $V_{_\Omega}(\Lambda)$. Then $x\in\Lambda\,$ iff
$\,\varphi_{_\Omega}(x)\odot\varphi_{_\Omega}(y_i)=0,\,\, \forall \,
i=1,\dots,m$.
\end{lemma}
\begin{proof}
Let $x \in \Lambda$. Since $\mathcal{B}\subseteq\Lambda$, by
orthogonality we have $\left\langle e_x,  e_{y_i}\right\rangle = 0,
\, \forall \,\,{y_i} \in \mathcal{B}$ and the result follows from
Lemma \ref{equivalence}.
\smallskip

For the converse, let $dim(V_{_\Omega}(\Lambda))=m$ and
$\mathcal{B}=\left\{y_1,\dots,y_m\right\}\subseteq \Lambda$ be such
that $\varphi_{_\Omega}(\mathcal{B})$ is a basis for
$V_{_\Omega}(\Lambda)$. Suppose there exists some $x\not\in\Lambda$
such that $\varphi_{_\Omega}(x)\odot\varphi_{_\Omega}(y_j)=0,
\forall y_j \in\mathcal{B}$.  Since $\varphi_{_\Omega}(\mathcal B)$
is a basis for $V_{_\Omega}(\Lambda)$, we have for any $\lambda \in
\Lambda$,  $\varphi_{_\Omega}(\lambda)=\sum_{j=1}^m a_j \
\varphi_{_\Omega}(y_j)$. Now by linearity we get
$\varphi_{_\Omega}(x)\odot\varphi_{_\Omega}(\lambda)=\sum_{j=1}^m
\overline{a_j}\
(\varphi_{_\Omega}(x)\odot\varphi_{_\Omega}(y_j))=0$. Hence by Lemma
\ref{equivalence} we get $\left\langle e_x,  e_\lambda\right\rangle
= 0, \,\,\forall \, \lambda \in \Lambda$. But
$E_\Lambda=\left\{e_\lambda: \lambda \in \Lambda \right\}$ is total
in $L^2(\Omega)$, and $e_x\not\equiv 0$, a contradiction. Thus $x$
must be in $\Lambda$.
\end{proof}

The following Lemma, gives a rather nice criterion for a spectrum
$\Lambda$ to be periodic.

\begin{lemma}\label{basis pattern}
Let $dim(V_{_\Omega}(\Lambda))=m\leq n$ and
$\mathcal{B}=\left\{y_1,\dots,y_m\right\}\subseteq\Lambda$ be such
that $\varphi_{_\Omega}(\mathcal{B})$ is a basis for
$V_{_\Omega}(\Lambda)$. If for some $d \in \R$, we have
$\mathcal{B}+d=\left\{y_1+d,\dots,y_m+d\right\}\subseteq\Lambda$
then $\Lambda$ is $d$-periodic, i.e., $\Lambda=
\{\lambda_1,\dots,\lambda_d\} +d \Z$.
\end{lemma}

\begin{proof}
By Lemma \ref{local finite} $x\in\Lambda$ iff
$\varphi_{_\Omega}(x)\odot\varphi_{_\Omega}(y_j)=0, \,j=1,\dots,m$.
Let $\lambda\in\Lambda$, since $\mathcal{B}+d\subseteq\Lambda$ we
get $\varphi_{_\Omega}(\lambda)\odot\varphi_{_\Omega}(y_j+d)=0,
\,j=1,\dots,m \iff
\varphi_{_\Omega}(\lambda-d)\odot\varphi_{_\Omega}(y_j)=0,
\,j=1,\dots,m \iff \lambda-d \in\Lambda$ and hence $\Lambda$ is
$d$-periodic. By Theorem \ref{D} we get $d\in\N$ and since $\Lambda$ 
has density $1$ by
Theorem \ref{landau}, we conclude that
$\Lambda= \{\lambda_1,\dots,\lambda_d\} +d \Z$.
\end{proof}

Recall, that if $\Gamma$ is periodic, has density $1$ and
$\varphi_{_\Omega}(\Gamma)$ is a set of mutually null-vectors, then
by \cite{P1}, \cite{LW2} $(\Omega,\Gamma)$ is a spectral pair.
\smallskip

Let $(\Omega,\Gamma)$ be such that $\varphi_{_\Omega}(\Gamma)$ is a
set of mutually null vectors. The natural queston is: Can we extend
$\Gamma$ to a spectrum of $\Omega$ ? The following Theorem gives a
criterion for periodic orthogonal extension of a set $\Gamma$  and
will be  central to our proof of periodicity of a spectrum in the
next section.

\begin{theorem}\label{orthogonal extension}
Let $\Gamma \subset \R$ be a such that the set of exponentials
$E_\Gamma$ are orthogonal in $L^2(\Omega)$. Let $dim
(V_{_\Omega}(\Gamma))=r$ and $\mathcal B_0=\{\mu_1,\dots,\mu_r\}$ be
such that $\varphi_{_\Omega}(\mathcal B_0)$ forms a basis of
$V_{_\Omega}(\Gamma)$. Further suppose a translate of $\mathcal B_0$
is contained in $\Gamma$, i.e., $\mathcal B_1=\mathcal B_0+d
\subseteq \Gamma$. Then $\Gamma$ can be extended periodically to
obtain a $d$-periodic subset $\Gamma_d\subseteq \R$ such that the
set of exponentials $E_{\Gamma_d}$ are orthogonal in $L^2(\Omega)$.
\end{theorem}

\begin{proof} Let $\Gamma_d:= \Gamma+ d\Z$. As in Lemma
\ref{basis pattern}, we will prove that
$\varphi_{_\Omega}(\Gamma_d)$ is a mutually null set. We will first
show by induction that
$$\varphi_{_\Omega}(\mu_k) \odot \varphi_{_\Omega}(\mu_j+ld) = 0
\mbox{ for all } l \in \Z, \mbox{ and } j,k=1,\dots,r.$$ Observe
that both $\varphi_{_\Omega}(\mathcal B_0)$ and
$\varphi_{_\Omega}(\mathcal B_1)$ span the same vector space
$V_{_\Omega}(\Gamma)$. Let us assume  that the
orthogonality relations hold for all $s=1,\dots,l-1$ i.e.,
$$\varphi_{_\Omega}(\mu_k)\odot\varphi_{_\Omega}(\mu_j+sd)=0
\mbox{ for all } j,k=1,\dots,r.$$
We have to show
$$\varphi_{_\Omega}(\mu_k)\odot\varphi_{_\Omega}(\mu_j+ld)=0 \mbox{
for all } j,k=1,\dots,r.$$ But by the induction hypothesis, we have
$$\varphi_{_\Omega}(\mu_k+d)\odot\varphi_{_\Omega}(\mu_j+ld)=
\varphi_{_\Omega}(\mu_k)\odot\varphi_{_\Omega}(\mu_j+(l-1)d)=0,\,
\forall j,k=1,\dots,r$$
But, we know that $\varphi_{_\Omega}(\mathcal B_0) \subseteq span
\{\varphi_{_\Omega}(\mathcal B_1) \}$. Hence,
$$\varphi_{_\Omega}(\mu_k)\odot\varphi_{_\Omega}(\mu_j+ld)=0,\,
\forall j,k=1,\dots,r$$

Now if $\gamma,\gamma'\in \Gamma_d$, then $\gamma=\gamma_p+ld$,
$\gamma'=\gamma_p'+l'd$ for some $\gamma_p,\gamma_p' \in \Gamma$ and
$l, l'\in\Z$. Since $\varphi_{_\Omega}(\gamma_p),
\varphi_{_\Omega}(\gamma_p') \in V_{_\Omega}(\Gamma) =
Span\left\{\varphi_{_\Omega}(\mathcal B_0)\right\}$, we have

$$ \varphi_{_\Omega}(\gamma_p)=\sum_{k=1}^{r} \alpha_k
\varphi_{_\Omega}(\mu_k) \mbox{ and } \varphi_{_\Omega}
(\gamma_p')=\sum_{j=1}^{r} \alpha_j' \varphi_{_\Omega}(\mu_j)$$

Now,

$$\varphi_{_\Omega}(\gamma)\odot\varphi_{_\Omega}(\gamma')=
\varphi_{_\Omega}(\gamma_p+ld) \odot \varphi_{_\Omega}
(\gamma_p'+l'd)=\varphi_{_\Omega}(\gamma_p+(l-l')d)\odot
\varphi_{_\Omega}(\gamma_p')$$

$$=\varphi_{_\Omega}(\gamma_p+(l-l')d) \odot \left
(\sum_1^{r}\alpha_j'\varphi_{_\Omega}(\mu_j)\right)$$

$$=\sum_{j=1}^{r} \overline{\alpha_j'} \varphi_{_\Omega}
(\gamma_p) \odot \varphi_{_\Omega}(\mu_j+(l'-l)d)$$

$$=\sum_{j=1}^{r} \sum_{k=1}^{r} \overline{\alpha_j'}
\alpha_k \varphi_{_\Omega}(\mu_k) \odot \varphi_{_\Omega}
(\mu_j+(l'-l)d) = 0$$

\end{proof}
\medskip

\begin{remark}
Under the assumption of Lemma \ref{basis pattern} the $\Lambda_d$
obtained in Theorem \ref{orthogonal extension} is $\Lambda$ itself.
\end{remark}

\subsection{Density of the spectrum} Let $\Gamma\subset\R$ be a uniformly discrete
set. Then we define $n^+(R),\,\, n^-(R)$ respectively, as the
largest and smallest number of elements of $\Gamma$ contained in any
interval of length $R$, i.e.,

$$ n^+(R) = \max\limits_{x \in \R} \#
\{ \Gamma \cap [x-R,x+R] \}$$ $$n^-(R) = \min\limits_{x \in \R} \#
\{ \Gamma \cap [x-R,x+R] \}.$$
\medskip

A uniformly discrete set $\Gamma$ is called a {\it set of sampling}
for $L^2(\Omega)$, if there exists a constant $K$ such that
$\left\|f\right\|_2^2\, \leq K \,\sum_{\lambda \in \Lambda}|
\hat{f}(\lambda)|^2,  \,\,\, \forall f \in L^2(\Omega)$, and $\Gamma$
is called a {\it set of interpolation} for $L^2(\Omega)$, if for
every square summable sequence $\{a_\gamma\}_{\gamma \in \Gamma}$,
there exists an $f \in L^2(\Omega)$ with $\hat{f}(\gamma) =
a_\gamma$, $\gamma \in \Gamma$.
\medskip

Clearly if $(\Omega,\Lambda)$ is a spectral pair, then $\Lambda$ is
both a set of sampling and a set of interpolation for $L^2(\Omega)$.
The following result of Landau, regarding sets of sampling and
interpolation gives an estimate on the numbers $n^+(R)$ and $n^-(R)$
for a spectrum $\Lambda$, when $\Omega$ is a union of a finite
number of intervals.

\begin{theorem} (Landau \cite{Landau})\label{landau} Let $\Omega$ be
a union of a finite number of intervals with  total measure $1$, and
$\Lambda$ a uniformly discrete set. Then\\

\begin{enumerate}
    \item If $\Lambda$ is a set of sampling for $L^2(\Omega)$,
    $$n^-(R) \geq R - A\log^+R -B $$
    \item If $\Lambda$ is a set of interpolation for $L^2(\Omega)$,
    $$n^+(R) \leq R - A\log^+R -B $$
where $A$ and $B$ are constants independent of $R$
\end{enumerate}

\end{theorem}
\medskip

It follows from Theorem \ref{landau} that $\Lambda$ has asymptotic
density 1, that is
$$\rho(\Lambda) := \lim_{R\rightarrow \infty} \frac{\# \left(
\Lambda \cap[-R+x,R+x]\right)}{2R} = 1, \mbox{ uniformly in}\, x
\in\R.$$
\medskip

{\section{Proof of periodicity of the spectrum}

Once again in this section $\Omega\subset \R$ is a union of finitely
many intervals, $\Omega = \cup_1^n [a_j,a_j+r_j],\, \sum_{j=1}^n r_j
=1$. We assume that $\Omega$ is spectral with a spectrum $\Lambda$.
We will continue to use the notations introduced in section
\ref{section 2}.
\smallskip

We begin with some definitions.
\smallskip

Let $\Lambda=\left\{\lambda_j\right\}_{j\in\Z}$ where
$\lambda_j<\lambda_{j+1}$ and $\lambda_0=0$. Recall that the
consecutive distance set of $\Lambda$, namely
$$\Lambda_s=\left\{\lambda_{j+1}-\lambda_j :j \in \Z\right\} $$
is finite. So we can view $\Lambda$ as an infinite word with a
finite alphabet $\Lambda_s=\left\{d_1,d_2,\dots,d_l\right\}$. For a
finite word $W=\left[d_{j_1},d_{j_2},\dots,d_{j_n}\right], d_{j_i}
\in \Lambda_s $ we write $length(W)=\sum_{i=1}^n d_{j_i}$.
\medskip

Suppose $dim(V_\Omega(\Lambda))=m \leq n$ and let $\{\mu_1, \mu_2,
..., \mu_m\}$ be such that $\{\varphi_{_\Omega}(\mu_j), \,
j=1,2,\dots,m \}$ is a basis for $V_\Omega(\Lambda)$.
\medskip

Choose $L_0$ such that  $\left\{\mu_1,\mu_2,\dots\mu_m\right\}
\subseteq\left[0,L_0 \right]$ and then for any $L \geq L_0$,
partition $\R$ as
$$\R = \cup_{k\in\Z} \left[kL,(k+1)L\right).$$

Let
$$\Lambda_k^L=\Lambda\cap\left[kL,(k+1)L\right),$$

\smallskip
Now, for each $ k \in \Z, \,\, \Lambda_k^L$, corresponds to a finite
word of length at most $L$, and there are only finitely many, say
$N_L$, words of length at most $L$. Let
$$V_k^L=Span\left\{\varphi_{_\Omega}(\lambda): \lambda \in
\Lambda_k^L\right\}.$$

Let us first consider the  {\bf special case} that for some large
enough $L$ we have
\begin{equation}\label{dim=m}
dim(V_k^L)=m \mbox{ for every } k \in \Z.
\end{equation}

In this case, each $\Lambda_k^L$ has a set of $m$ elements $\mathcal
B_k:=\left\{\mu^k_1, \mu^k_2, \dots,\mu^k_m\right\}$ such that
$\varphi_{\Omega}(\mathcal B_k):=\left\{\varphi_{_\Omega}(\mu^k_1),
\dots,\varphi_{_\Omega}(\mu^k_m) \right\}$ forms a basis of
$V_{_\Omega}(\Lambda)$. Also by the remarks above, at least two of
the words $\Lambda_{k_1}^L$ and $\Lambda_{k_2}^L$ must be the same.
Hence for some $d \in \R,\,\, \Lambda_{k_2}^L = \Lambda_{k_1}^L + d
$. In particular, there exists  $k_0$, such that  $\Lambda_{k_0}^L$
contains a set of elements
$\left\{\mu^{k_0}_1,\dots,\mu^{k_0}_m\right\}$ which form a basis
of $V_{_\Omega}(\Lambda)$ and also
$\left\{\mu^{k_0}_1,\dots,\mu^{k_0}_m\right\}+d \subseteq \Lambda$.
Thus the hypothesis of Lemma \ref{basis pattern} holds, and so
$\Lambda$ is $d$-periodic.
\smallskip

Observe that in the above argument, we do not require as much as
(\ref{dim=m}). It would be enough if $\left\{k: dim(V_k^L) = m
\right\}$ is an infinite set, or for that matter, has at least
$N_L+1$ elements. But once we conclude that $\Lambda$ is periodic,
it will follow that for some, possibly larger $L'$,\,(if $d$ is the
period $L'=3d$ will do) that $dim(V_k^{L'})=m$, $\forall k \in \Z$.
\medskip

For the general case, let $1\leq s\leq m,$ and $L>0$ and write
$$E_s^L = \left\{k: dim(V_k^L) \geq s \right\}$$

We have just seen that if for some $L>0$, $E_m^L=\Z$, then $\Lambda$
is periodic. Suppose this is not the case. Then we need the following lemma:
\smallskip

\begin{lemma}Let $m' \leq m$ be the largest integer such that there
exists an $L'>0$ so that $E_{m'}^{L'} =\Z$. Then $m'$ itself will
occur infinitely often in the set $\left\{dim
(V_k^{L'})\right\}_{k\in\Z} $.
\end{lemma}

{\it Proof:} First note that for $s=1$, we can choose $L'
> max\{d_j\}$, and then $E_1^{L'}=\Z$ so clearly $m' \geq 1$. If
$dim(V^{L'}_k)=m'$ only for finitely many $k$'s then we can take
$\tilde{L}$ large enough so that $dim(V^{L'}_k)=m'$ for precisely
one interval of the partition $\{[k\tilde{L}, (k+1)\tilde{L}]\} $.
Let $L^{''} = 2 \tilde{L}$, then observe that $E^{L^{''}}_{m'+1}=\Z$
, and this contradicts maximality of $m'$. (Without loss of
generality we may choose $L' \in \N$.)
\medskip

We will now prove {\bf Theorem \ref{theorem 1}}.

\begin{proof}

{\bf Step 1.} We will first prove that the spectrum $\Lambda$ can be
modified to a set $\Lambda_d$ which is $d$-periodic and is such that
$(\Omega, \Lambda_d)$ is a spectral pair. For this we use Landau's
density result to extract a ``patch'' from $\Lambda$ which has some
periodic structure and has a large enough density. Then we use Theorem
\ref{orthogonal extension} to show that a suitable
periodization of this patch is a spectrum.
\medskip

With $L'$ as above, let
$$\epsilon_{L'}=\frac{1}{2L'(N_{L'}+1)}.$$
Then choose $L^{*}>\frac{1}{2\epsilon_{L'}}=L'(N_{L'}+1)$ such that
 $n^-(L^*)/L^*>1-\epsilon_{L'}$.\\
In the case under consideration, we know that $E_{m'}^{L^*}=\Z$ and
also that the cardinality of the set
$\left\{p:dim(V^{L^*}_p)=m'\right\}$ is infinite. We choose and fix
one such $p$ such that $dim(V_p^{L^*})=m'$. By the choice of $L^*$,
the interval $\left[pL^*,(p+1)L^*\right)$ contains at least
$(N_{L'}+1)$ disjoint intervals of length $L'$. Now for $j=1,...,
N_L+1$ each of the $\Lambda_j^{L'} \subset
\left[pL^*,(p+1)L^*\right)$ has a word $W_j$ of length at most $L'$
associated with it. Further, observe that by the choice of $L'$, each
of these $\Lambda_j^{L'}$ contains at least $m'$ elements whose
image under $\varphi_{_\Omega}$ is a linearly independent set, and
that, by the choice of $p$, there can be at most $m'$ such elements.
Notice this implies $V^{L^*}_p=V^{L'}_j,\,\, j=1,\dots,N_L+1$.
\medskip

Hence by the pigeon hole principle, there exists $k_1$ and $k_2$
such that the words $\Lambda_{k_1}^{L'}$ and $\Lambda_{k_2}^{L'}$
are the same, and therefore $\Lambda_{k_2}^{L'} =
\Lambda_{k_1}^{L'}+d$ for some $d \in \R,$ where
$d\leq(N_{L'}+1)L'=\frac{1}{2\epsilon_{L'}}$.

\medskip

To complete the proof, we will need the following lemma:

\begin{lemma}\label{2em}
Let $\Lambda_d$ be the $d$-periodization of $\Lambda_p^{L^*}$, i.e.
$ \Lambda_d =\left\{\Lambda_p^{L^*}+d\Z\right\}$. Then
 $\Lambda_d$ is orthogonal.
\end{lemma}

\begin{proof}  Let $\mathcal{B}_0=\left\{\mu_1,\dots,\mu_{m'}\right\}
\subseteq \Lambda_{k_1}^{L'} \subseteq \Lambda_p^{L^*}$ be such that
 $\varphi_{_\Omega}(\mathcal B_0):=\{\varphi_{_\Omega}(\mu_1),\dots,
 \varphi_{_\Omega}(\mu_{m'})\}$
is a basis of $V_p^{L^*}$ and also of $V_{k_1}^{L'}$. Now since
$\Lambda_{k_1}^{L'}+d = \Lambda_{k_2}^{L'}, \,\,\mathcal{B}_1=
\mathcal{B}_0+d \subseteq \Lambda_{k_2}^{L'}$,  this subset  again
gives a basis for $V^{L^*}_p$. By Theorem \ref{orthogonal extension}
we see that the set of exponentials $E_{\Lambda_d}$ are mutually
orthogonal in $L^2(\Omega)$.
\end{proof}
 Now since $\Lambda_d$ is orthogonal it is a set of interpolation
 and by Landau's density theorem we get $\rho(\Lambda_d)\leq 1$.
 But by our choice of $L^*$ we get $\rho(\Lambda_d) > n^-(L^*)/L^* >
 1- \epsilon_{L'}$. On the other hand, since $\Lambda_d$ is d-periodic,
 if $\rho(\Lambda_d)<1$, we have $\rho(\Lambda_d)\leq 1-\frac{1}{d}
 <1-2\epsilon_{L'}$ as $\frac{1}{d}\geq 2\epsilon_{L'}$. This is a
 contradiction.
\smallskip

It follows that $\Lambda_d$ is a periodic set whose density is $1$ and
$E_{\Lambda_d}$ is orthogonal in $L^2(\Omega)$. Thus we get
$\Lambda_d$ is a spectrum for $\Omega$ \cite{P1},\cite{LW2}. Since
$\Lambda_d$  has density of $1$ and is $d$-periodic it can be
written in the form $\Lambda_d = \cup_{j=1}^d (\mu_j +d \Z)$.
\medskip

{\bf Step 2.} We now prove that $\Lambda$ itself is periodic. Once
again we will be using Landau's density Theorem and Theorem
\ref{orthogonal extension} along with Theorem \ref{D} which will be
crucial.
\smallskip

Choose $L^*$ as above, so that $\left\{p:dim(V^{L^*}_p)=m'\right
\}$ is infinite.\\

Then let $L^{**}$ be such that $$n^-(L^{**})/L^{**} > 1-\frac{1}
{2(n+1)L^*} \mbox{ and }L^{**}>>(n+1)L^*$$ (Recall that $n$ is the
 number of intervals in $\Omega$ ). Here by $>>$ we mean that many
 blocks of intervals, each of length $(n+1)L^*$ are contained in
 any interval of the $L^{**}-$grid.
\medskip

Then we can find a $p$ such that $dim(V^{L^{**}}_p)=m'$ (since
there are infinitely many such). Now extend   $\Lambda_p^{L^{**}}$
$ d-$periodically to a spectrum $\Lambda_d^*$ of $\Omega$, where
$d<L^*$. Write $\Lambda_d^* = \bigcup_{j=1}^{d} (\mu_j + d\Z)$,
with $\mu_1, \mu_2,..., \mu_d \in \left[pL^{**},(p+1)L^{**}\right)$.

\medskip

We end the proof by showing that in fact $\Lambda_d^* = \Lambda$.
For this it will be enough to show that for each $\mu_j$, there are
$(n+1)$ consecutive terms from the arithmetic progression $\mu_j
+d\Z$ in $\Lambda$. Suppose this is not the case, then for each
$a\in\Z$ such that $[\mu_j+a d,\mu_j+(a+n)d] \subset [p L^{**},(p+1)
L^{**}]$ we have at least one element from the n+1 length AP
$\mu_j+a d,\mu_j+(a+1)d\dots,\mu_j+(a+n)d$ is missing from
$\Lambda_p^{L^{**}}$. But that will affect the density, so that
$n^-(L^{**})/L^{**} \leq 1-\frac{1}{(n+1) d} \leq
1-\frac{1}{(n+1)L^*}$, which is a contradiction. Now By Theorem
\ref{D} we get that $\Lambda$ is indeed periodic.
\end{proof}

The structure of spectral sets $\Omega$ which have a periodic
spectrum is well known (see \cite{P1},\cite{LW2}). Here for the sake
of completeness we give a structure theorem for $\Omega$ using a
result of Kolountzakis.

\medskip
{\bf Theorem (Kolountzakis)}.{\it Let $\Omega$ be a bounded open
set, $\Lambda$ a discrete set in $\R^d$,  and $\delta_\Lambda
=\sum_{\lambda \in \Lambda} \delta_\lambda$. Then $|\widehat{
\chi_\Omega}|^2 + \Lambda $ is a tiling if and only if $\Lambda$ has
uniformly bounded density and $$(\Omega - \Omega) \cap  supp(
\widehat{\delta_\Lambda}) =\{0\}.$$}

We will now prove {\bf Theorem \ref{theorem 2}}.

\begin{proof}

Recall that
$(\Omega,\Lambda)$ is a spectral pair if and only if $|\widehat{
\chi_\Omega}|^2 + \Lambda $ is a tiling. Further if $\Lambda$ is
$d$-periodic, then $\Omega$ $d$-tiles $\R$, i.e.
$\sum_n\chi_\Omega(x+n/d)=d.$ (i.e. the set $\Z/d$ $d$-tiles $\R$ by
$\Omega$.)

In particular,
$$ d \chi_{[0,\frac{1}{d})}(x) = \chi_{[0,\frac{1}{d})} \sum_{k \in \Z}
\chi_{\Omega} (x + k/d).$$

So for each $x \in [0,\frac{k}{d})$, the set $A_x=\{k \in \Z :
x + k/d \in \Omega\}$ has cardinality $d$.
Define an equivalence relation $\approx$ on $[0,1/d)$ by
$x\approx y$ if and only if $A_x = A_y$.

Since $\Omega$ is bounded, the above equivalence relation gives a
partition of $[0,1/d)$ into finitely many equivalence classes
$E_1,E_2,\dots,E_k$. For each $E_j$ and we write $A_j$ for the
common set defined above.
\medskip

Then $\Omega = \cup_{j=1}^k (E_j + A_j)$ and $[0,1/d)= \cup_{j=1}^k
E_j$ and we may assume $|E_j|>0 \ \ \forall j = 1,2,\dots,k$. Now
let $\Omega_j := [0,1/d) + A_j$. Our claim is $(\Omega_j,\Lambda)$
is a spectral pair. We will need the following theorem due to
Kolountzakis \cite{K1}.

Now as $\sharp(A_j)=d$ we have
$|\Omega_j|=1$. If $\Lambda = \Gamma + d \Z$ with $\Gamma =
\{\lambda_1, \lambda_2, ..., \lambda_d\}$, then
$supp(\widehat{\delta_\Lambda}) = \{ k/d : \widehat{\delta_\Lambda}
(k/d) \neq 0 \} \subseteq \Z/d $ and $supp(\Omega' - \Omega')
\subseteq (- 1/d, 1/d ) + A_j - A_j $. But $A_j - A_J$ is
$1/d$-separated, so $supp(\Omega' - \Omega') \cap
supp(\widehat{\delta_\Lambda}) = \{0\}$ for otherwise as $E_j+A_j
\subseteq \Omega$ and $|E_j|>0$ we get $supp(\Omega - \Omega) \cap
supp(\widehat{\delta_\Lambda}) \neq \{0\}$ and thus
$(\Omega,\Lambda)$ cannot be a spectral set.
\end{proof}

\bigskip
{\bf Acknowledgement.} The authors would like to thank Krishnan Rajkumar and C.P. Anil Kumar for the many insightful comments and suggestions they made at several stages of this work and for providing us with much needed encouragement.

\end{document}